\documentclass[journal,twoside,web]{ieeecolor}
\usepackage{generic}
\usepackage{cite}
\usepackage{amsmath,amssymb,amsfonts}
\usepackage{algorithmic}
\usepackage{graphicx}
\usepackage{textcomp}
\def\BibTeX{{\rm B\kern-.05em{\sc i\kern-.025em b}\kern-.08em
    T\kern-.1667em\lower.7ex\hbox{E}\kern-.125emX}}
\usepackage{graphicx}
\usepackage{psfrag}
\usepackage{centernot}
\usepackage{units}
\usepackage{nicefrac}
\usepackage{multirow}
\usepackage{fancyhdr}
\usepackage{rotating}
\usepackage{multirow}
\usepackage[auth-sc,affil-sl]{authblk}
\usepackage{amsmath,amssymb}   
\usepackage{inputenc} 
\usepackage[english]{babel}    
\usepackage{caption}     
\makeatletter
    \let\NAT@parse\undefined
    \makeatother
\usepackage[square,numbers]{natbib}     
\renewcommand\thesection{\arabic{section}}

\newcommand{\ie}{\emph{i.e.\ }}
\newcommand{\eg}{\emph{e.g.\ }}
\newcommand{\equalperdefinition}{\ensuremath{\stackrel{\text{\tiny def}}{=}}}
\newcommand{\matrixminus}{\ensuremath{\textit{--}}}
\newcommand{\matrixminusspace}{\ensuremath{\phantom{\textit{{--}}}}}
\newcommand{\bs}[1]{\ensuremath{\boldsymbol{#1}}}

\DeclareMathOperator*{\linspan}{\ensuremath{\mathrm{span}}}

\newtheorem{algo}{Algorithm}
\allowdisplaybreaks

\begin{document}
\title{Automatic traction control for articulated off-road vehicles} 
\author{Johan Markdahl
\thanks{This paragraph of the first footnote will contain the date on 
which you submitted your paper for review. This work was jointly supported by the Swedish innovation agency VINNOVA and Volvo Construction Equipment.}
\thanks{Johan Markdahl is with the Luxembourg Centre for Systems Biomedicine at the University of Luxembourg (e-mail: markdahl@kth.se).}}

\maketitle

\begin{abstract}
Construction equipment is designed to maintain good traction, even when operating in difficult off-road conditions. To curb wheel slip, the vehicles are equipped with differential locks. A driver may engage/disengage the locks to switch between two distinct operating modes: the closed mode is characterized by greater off-road 
passability while the open mode allows better manueverability. However, many drivers lack the education and experience required to correctly judge the terrain ahead of the vehicle and therefore engage/disengage the locks in a suboptimal fashion. An automatic traction control solution for locking and opening the differentials is hence desirable. This paper compares three on/off differential lock control algorithms, all derived from the same kinematic vehicle model but each relying on the availability of output signals from different sensors. The validity of the kinematic model and the algorithms' sensitivity to the values assumed by a couple of unobservable states, the wheel slip angles, is investigated by comparison to a realistic articulated hauler model in the multibody physics simulator MSC ADAMS.
\end{abstract}


\section{Introduction}

\noindent The design of articulated haulers is optimized for traversing difficult terrain \cite{terrangmaskinen,history}. This allows the vehicles to take the shortest route on a load, haul, dump run; thereby minimizing fuel consumption and time expenditure. Part of this optimization is automatic traction control (ATC) which eases the decision making burden of the driver, protects the tires from unnecessary wear, and reduces fuel consumption by up to 6\% \cite{atc}. Note that the pricetag for a single tire is USD \mbox{5 000} and tires represent 20--25\% of a haulers operating costs \cite{ulf}. Manufacturers like Caterpillar, John Deree, Komatsu, and Volvo  use ATC \cite{article}, but the details of their algorithms are unknown to the public except for glimpses gained from ads \cite{atc} and patents \cite{hosseini1996method,smith2000differential,holt2003electronic,murakami2003control,olsson2008,smith2010differential,newberry2007differential}. The ATC problem is also interesting from a theoretical point of view since the nature of the actuators (so-called dog clutches) requires on/off control with a strong emphasis on the `on' decision.  This paper proposes three traction control algorithms based on different sensor output. The algorithms are validated against a realistic model of a Volvo Construction Equipment A40X articulated hauler in the multibody physic simulator MSC ADAMS \cite{illerhag2000study}.


\subsection{Background}
\label{sec:prob}

\noindent Articulated haulers are heavy equipment, used to transport large quantities of loose materials such as sand, gravel, and liquids in off-road enviroments. The articulated steering (see Fig \ref{fig:haulerBirdsview}) ensures high maneuverability although at the cost of a lower maximum payload. Articulated haulers may, for example, be employed at mines to transport ore over grounds not traversable by ordinary vehicles or on construction sites to transport building material.

Traction is an adhesive friction force in tire/road interface that serves to drive the vehicle forward. A tire may sometimes lose its grip and slip rather than roll over the road, \eg\ if subject to full throttle, icy conditions or a steep inclination. This is referred to as lost traction. Lost traction is undesirable since it reduces the vehicle's traversability and increases tire wear. Means must therefore be taken to curb wheel slip and regain lost traction, preferably at the onset of wheel slip.

A differential is a driveline component that distributes power, \ie torque and rotational speed, from an input shaft to two output shafts. The differentials on articulated haulers have two distinct operating modes: open and closed. An open differential distributes rotational speed freely and torque evenly. A locked differential forces the output shafts to assume the lowest of the two wheel speeds while torque is distributed freely. A succesion of locks may be engaged to curb the slip of multiple wheels. 

The differential locks are of the dog-clutch variety: a pair of face gears that are locked together pneumatically and pulled apart by a spring. Being locked together, the two shafts assume the same speed. The dog-clutch has a range of angular velocity differences over which it is safe to engage, otherwise engagement risks damaging the gear teeth. If the angular speed differences does not satisfy this constraint, then some articulated hauler models allow individual wheel brakes to reduce the rotational speed of selected shafts.


The articulated hauler often serves as an entry point for beginner drivers who later move on to more advanced vehicles such as wheel loaders and excavators as they gain experience. Inexperienced drivers tend to be overtly reliant on the differential locks, turning them on at all times. By contrast, system logs from Volvo's articulated haulers show that their ATC system can outperform even skilled drivers.

\subsection{Problem statement}

\noindent The key question for ATC is when to lock the differentials. This comes down to comparing the revolutions per minute (RPM) of driveline shafts that are coupled through open differentials. A difference in RPM indicates slip, unless the vehicle is turning or braking. The most interesting automatic control challange is to discern differences in RPM due to turning from those due to wheel slip. In this paper we use the hauler geometry to develop a kinematic model that lets us answer this question. Moreover, we explore the advantages that may be gained from utilizing information obtained from a ground speed sensor (\eg a ground speed radar or a GPS receiver) and individual wheel tachometers (angular speed sensors). These sensors are not standard, so this paper also addresses questions that are of interest to the manufacturers.

\subsection{Literature review and contribution}

\noindent The literature on traction control and various related topics such as antilock breaking systems, electronic stability control, and anti-slip regulation mainly concerns automobiles, see \eg the survey \cite{pretagostini2020survey}, and the survey \cite{ivanov2014survey} on traction control for electric vehicles. Traction control for articulated off-road vehicles is largely unexplored in the academic literature, although there are some  exceptions \cite{ulf}. Rather, the state-of-the-art exists as inhouse software solutions that are unavailable to the public. To gain an idea of what ATC algorithms are used by manufacturers we may turn to commercials, product specifications, white papers, grey literature, and patents.

Some algorithms presented in patents are based on braking the slipping wheel, including \cite{hosseini1996method,holt2003electronic}. This actuation is not without limitations, including additional wear on the brakes. It requires sensors that can detect which in a pair of wheels that is slipping, \ie individual wheel tachometers. Moreover, the heat generated by braking can only be sustained by the vehicle from a limited time \cite{newberry2007differential}. The majority of relevant patents use differential locks for actuation \cite{smith2000differential,murakami2003control,newberry2007differential}. The algorithm \cite{murakami2003control} is unique in that it uses detection of oscillations in the driveline shafts to detect the presence or absence of wheel slip. 

Most algorithms make use of a steering angle sensor  \cite{hosseini1996method,holt2003electronic,murakami2003control,newberry2007differential,olsson2008} or a related output signal \cite{smith2000differential}, but unlike our paper they do not use the steering angle derivative. For motor graders an additional angle sensor for the front wheel pair is needed \cite{newberry2007differential,smith2010differential}. The algorithm \cite{smith2000differential} uses data from presure sensor in the steering hydraulics to decide when to lock or unlock the differentials. This algorithm also assumes that readings from a  ground speed sensor are  available. Some patents use wheel speed sensors \cite{hosseini1996method,holt2003electronic}, however, that does not necessarily mean that such algorithms are employed at present.

ATC by means of locking differentials have not received much attention in academic literature. The contribution of this paper is to provide three novel traction control algorithms based on a kinematic model of an articulated hauler. The kinematic model and the algorithms effectiveness are proved by comparison to a realistic articuled hauler model in MSC ADAMS. The results also apply to wheel loaders, which are obtained as a special case where the geometry satisfies $l_1=l_2$, see Section \ref{sec:hauler} and Figure \ref{fig:haulerBirdsview}. Some preliminary results from this paper appear in the author's master's thesis \cite{markdahl2010}. A draft version of this paper also appears in the PhD thesis \cite{ulf}.

\section{Preliminaries}

\subsection{Articulated hauler model in ADAMS}

\label{sec:hauler}

\noindent The results of this study are validated by simulations of an articulated hauler model in MSC ADAMS which is owned by Volvo Construction Equipment (VCE). Part of the presentation is therefore based on VCE articulated haulers (note that the author decleares no conflict of interests). In particular, assumptions about the vehicle geometry, the sensor network and actuators that are available on all units are based on the VCE ATC system \cite{atc}. However, the key ideas of our control algorithms generalize to other geometries and sensor networks.

\begin{figure}[htb!]
	\centering
	\includegraphics[width=0.50\textwidth, clip, trim=0cm 0cm 0 0.5cm]{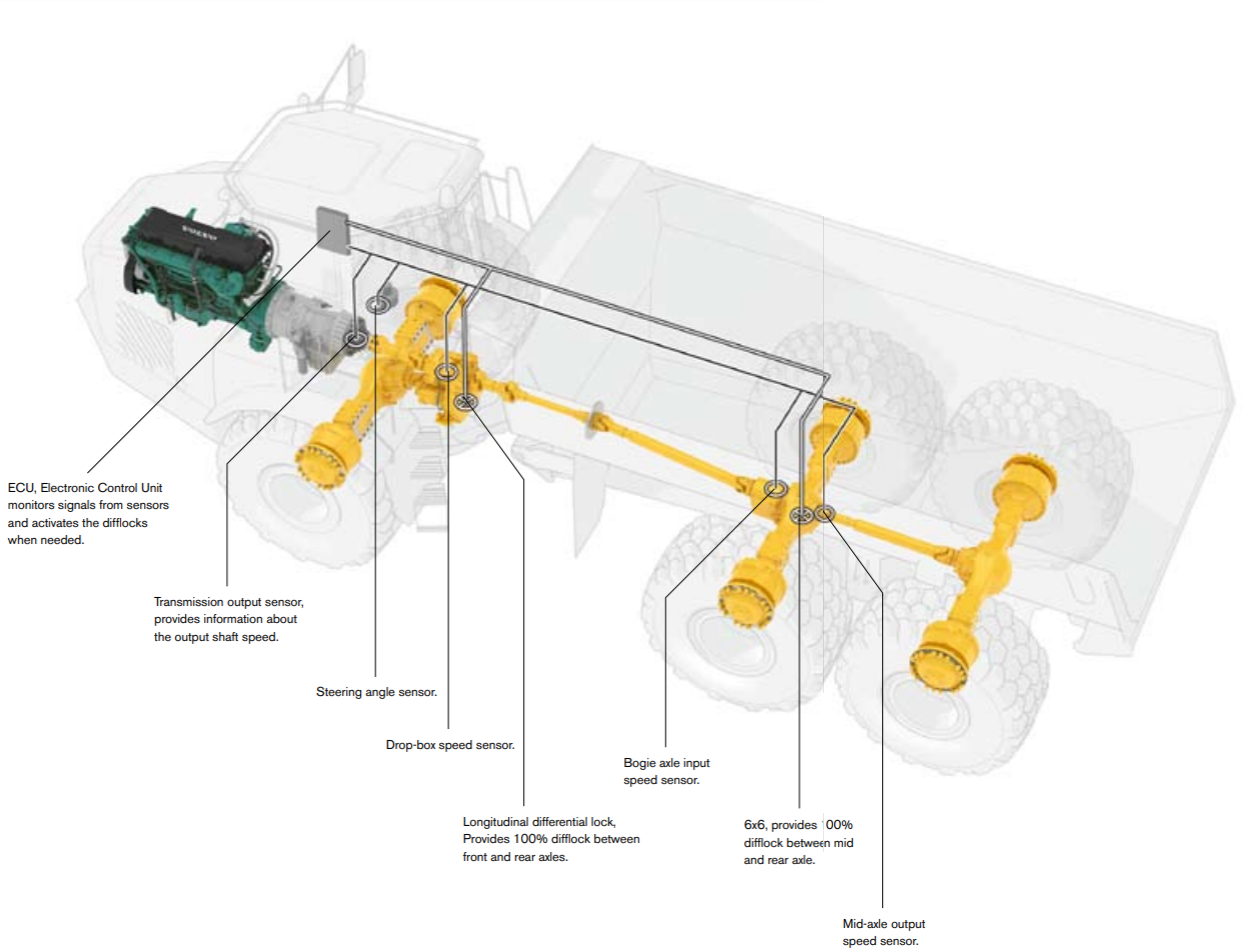} 
	\caption[ATC system sensor configuration.]{VCE ATC system sensor and actuator network \cite{atc}.}
	\label{fig:atc}
\end{figure}

The VCE sensor network is displayed in Fig. \ref{fig:atc}. The following five sensors are availible on all hauler units: a steering angle sensor; four tachometers: the transfer case (dropbox) in/out sensors and the bogie axle in/out sensors. The following five actuators are available: three transversal differential locks, one for each wheel axle (not in Fig. \ref{fig:atc}); one longitudinal differential lock which locks together the tractor and trailer unit axles; and one bogie lock which engages the $6\times6$ drive (the vehicle usually drives on $6\times4$ for better fuel economy).

\begin{figure}[htb!]
	\centering
	\includegraphics[width=0.50\textwidth, clip, trim=0cm 0cm 0 0.5cm]{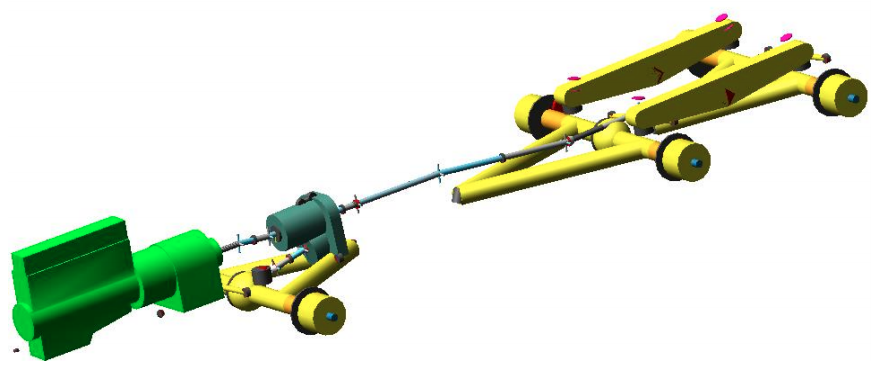} 
	\caption{The powertrain in the ADAMS model \cite{furmanik2015}.}
	\label{fig:adams}
\end{figure}

The geometry of the MSC ADAMS vehicle is the same as the A40X Volvo CE articulated hauler (X denotes the generation, which varies from D to G). Some relevant specifications:
\begin{itemize}
	\item Vehicle mass is 28 500 kg.
	\item Max load is 39 000 kg.
	\item Front wheel axle to steering joint, $l_1=1.278$ m. Front boogie wheel axle to steering join, $l_2=3.265$ m.
	\item The axle track, \ie the distance between a pair of wheels, is $2c_1=2.636$ m \cite{spec}.
\end{itemize}
A picture of the ADAMS model powertrain corresponding to the VCE A40X hauler is displayed in Fig. \ref{fig:adams}.

\subsection{Wheel slip definitions}

\noindent The (longitudinal) wheel slip $\lambda_l$ is often defined as a normalized function of the difference between the speed of the tire circumference and the ground, \emph{e.g.,}
\begin{align}\label{eq:slip}
\lambda_l\equalperdefinition \begin{cases}\frac{\omega r-v\cos\alpha}{\omega r} &\textit{ if }\omega r\geq v\cos\alpha\\
\frac{v\cos\alpha-\omega r}{v\cos\alpha} &\textit{ if }\omega r\leq v\cos\alpha
\end{cases}
\end{align}
where $\omega$ is the wheel angular speed, $r$ the outer tire radius, $\alpha$ the tire slip angle, and $v$ the ground speed, see \cite{terrangmaskinen,nielsen,andreev2010driveline}. The two cases in \eqref{eq:slip} correspond to driving and braking respectively. Together they imply $\lambda_l\in[0,1]$. We get back to wheel slip angles $\alpha$ in Section \ref{sec:general}, see also Fig. \ref{fig:axle} and \ref{fig:haulerBirdsview}. For now we just consider longitudinal slip. For an alternative, simpler longitudinal slip definition, consider the non-normalized quantity
\begin{align}
s_l\equalperdefinition\omega r-v\cos\alpha, \label{eq:deltaomega}
\end{align}
where $s_l\in\mathbb{R}$ can be interpreted as the distance a slipping wheel slips per second along its direction of orientation. 

The normalized slip quantity $\lambda$ is of interest since it relates to the tire/road interface friction coefficient $\mu(\lambda)$, as detailed in a body of empirical studies, see \emph{e.g.,} the survey \cite{li2006integrated}. The friction coefficient $\mu(\lambda)$ in turn enters vechicle dynamics including simple models such as the one-wheel and bicycle models as well as more advanced models \cite{nielsen}. Since the algorithms developed in this paper does not use tire dynamics, there is little benefit in adopting $\lambda$  as a measure of slip.

The non-normalized slip quantity \eqref{eq:deltaomega} has the advantage of being more directly related to tire wear. In fact, some amount of wheel slip always exists in the tire/road interface, but small to moderate amounts of wheel slip are tolerable. To see how small and large amounts of slip affects our definitions let $\alpha=0$, $r=1$ and consider two situations: 
\begin{itemize}
	\item[(i)] $\omega=\varepsilon$, $v=0$ and
	\item[(ii)] $\omega=2M$, $v=M$, 
\end{itemize}
where $\varepsilon \ll M$. Note that (ii) may result in significant tire wear whereas small amounts of slip as in case (i) are tolerable. Definition \eqref{eq:slip} and $r=1,\,\alpha=0$ yields $\lambda_l=1$ in case (i) and $\lambda_l=1/2$ in case (ii). Definition \eqref{eq:deltaomega} yields a $s_l=\varepsilon$ in case (i) and a $s_l=M$ in case (ii). The definition \eqref{eq:deltaomega} is preferable since it clearly distinguishes between the  case (i) and (ii) in a way that captures the fact that (ii) could potentially result in significant tire wear whereas (i) could not.

\newpage
\section{Main results}

\subsection{Control strategy}
\label{sec:control}

\noindent Certain sensor output signals are barely affected by wheel slip levels while others may change rapidly. Ground speed calculated from a GPS receiver's position readings is an example of an unaffected output. Tachometer measurements of the rotational speed of a drive shaft is an example of an affected output. This notion of either a static or a transient behavior of sensor output signals at the onset of wheel slip forms the basis of the traction control algorithms in this paper.

The rotational speed $\omega$ of a driveline shaft and the ground speed $v$ are proportional to each other under the assumption of zero slip and a rigid driveline:
\begin{align}
\omega r/i=v, \label{eq:proportional}
\end{align}
where $i$ is a gear conversion ratio (possibly equal to 1) and $r$ the outer tire radius. In the presence of wheel slip, the relation \eqref{eq:proportional} is replaced by our definition \eqref{eq:deltaomega} of wheel slip, $s=\omega r/i-v\cos\alpha$. We can generalize this idea to other kinematic equations that involve the vehicle velocity.

Let $\bs{y}\in\mathbb{R}^k$ be the sensor output signals and let $\bs{z}\in\mathbb{R}^l$ denote relevant states that are not measured. Let $g(\bs{y},\bs{z}):\mathbb{R}^k\times
\mathbb{R}^l\rightarrow\mathbb{R}$ be a function of the states and suppose that a kinematic equation
\begin{align}
g(\bs{y},\bs{z})|_{\bs{s}=\bs{0}}=0, \label{eq:equation}
\end{align}
holds in the absence of slip. The equality \eqref{eq:equation} is not guaranteed to withhold while $\bs{s}\neq\bs{0}$. Any difference in the left- and right-hand side of equation \eqref{eq:equation} can therfore be used as an indicator of wheel slip. Our approach to traction control is hence: if the left- and right-hand side of equation \eqref{eq:equation} differ beyond some preset tolerance, then engage a differential lock. 

\begin{figure}[hbt!]
	\centering
	\includegraphics[width=0.4\textwidth]{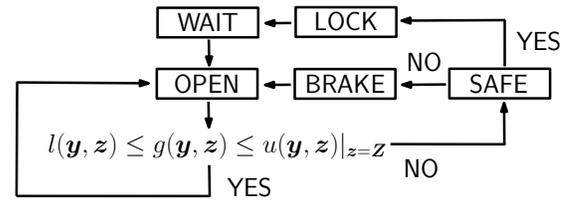}
	\caption{\label{fig:strategy}Schematic control strategy. The focus of this paper is on the condition $l(\bs{y},\bs{z})\leq g(\bs{y},\bs{z})\leq u(\bs{y},\bs{z})|_{\bs{z}=\bs{Z}}$ for locking.} 
\end{figure}

The basic idea of our control strategy is summarized in Fig. \ref{fig:strategy}. The controller monitors the value of $g(\bs{y},\bs{z})$ to see if slip beyond a preset treshhold is detected. Before the locks are engaged, a check is performed to ensure that the difference in angular velocity over the dog clutch is not too large. Engaging the clutch while the difference is large risks damaging the teeth. Hence apply the brakes prior to locking. After locking, it may not be possible to disengage the locks due to wind-up torques in the driveline. The dog clutch remains locked until the conditions are such that the locks can be disengaged.

Equation \eqref{eq:equation} is likely to be a kinematic relation. Most manufacturers own patented traction control algorithms based on articulated vehicle kinematics \cite{hosseini1996method,smith2000differential,holt2003electronic,murakami2003control,olsson2008,smith2010differential}. Ideally we would be able to measure the rotational and ground speed of each wheel, setting $\bs{y}=[\omega_i,v_i]^\mathsf{T}$ and $g_i(\bs{y})=\omega_ir_i-v_i, i\in\{1,\ldots,6\}$. However, neither individual wheel tachometers nor ground speed sensors form a part of the basic sensors configuration on VCE articulated haulers. As such we need to consider alternatives.

As a more feasible example, consider
\begin{align*}
\omega_{\textit{dbx,\,in}}-\omega_{\textit{dbx,\,out}}=\left\{ \begin{aligned}
&0 &&\mathrm{if} \, \gamma,s_i=0,  \forall i\in\{1,\ldots,6\}\\
&g(\bs{y},\bs{z}) &&\mathrm{otherwise,}
\end{aligned} \right.
\end{align*}
where $\omega_{\textit{dbx, in}}$ and $\omega_{\textit{dbx, out}}$ are the dropbox in/out shafts angular speeds (see Figure \ref{fig:atc}), $\gamma$ is the steering angle between the tractor and trailer unit, and $g$ is some function. However, note that if $\gamma=0$ and a tractor and trailer wheel should slip simultaneously, the above equation could still hold. This illustrates a limitation in our control design strategy: $\bs{s}=0\Rightarrow g(\bs{y},\bs{z})=0$ but $g(\bs{y},\bs{z})=0\centernot{\Rightarrow}\bs{s}=0$.

Introduce the functions $u(\bs{y})$ and $l(\bs{y})$ and relax the no slip constraint \eqref{eq:equation} as
\begin{align}
\begin{aligned}
g(\bs{y},\bs{Z})|_{\bs{s}=0}&\geq l(\bs{y}), \\
g(\bs{y},\bs{Z})|_{\bs{s}=0}&\leq u(\bs{y}),
\end{aligned} \label{eq:condition}
\end{align}
where $\bs{Z}$ is a constant nominal value of the unknown quantity $\bs{z}$. The functions $l$ and $u$ are included to account for the errors arising from setting $\bs{z}=\bs{Z}$. Note that there is a trade-off between functions $l,u$ with small magnitudes, which make control action fast but increase the risk of unnecessary engagements and functions $l,u$ with larger magnitudes, which reduce the risk of unnecessary engagements but also delay control action. Ideally, $g(\bs{y},\bs{z})=g(\bs{y})$ so that $u$ and $l$ may be set to small values. If the chosen expressions $g,l,u$ depends on unobservable states $\bs{z}$, then setting $\bs{z}=\bs{Z}$ may require the functions $u,l$ to have large magnitudes. In practice, as we see in Section \ref{sec:tuning}, $u,l$ can be defined as piece-wise linear functions of $\gamma$ based on data from simulations.

Finally, we provide a formal statement that encompasses the family of traction control algorithms we consider:

\begin{algo}
	Let	$g(\bs{y},\bs{z})|_{\bs{s}=0}=0$ be a kinematic relation that holds in the absence of slip. Suppose that we have two bounds $l(\bs{y},\bs{z})$ and $u(\bs{y},\bs{z})$ such that
	\begin{align*}
	l(\bs{y},\bs{z})|_{\bs{s}\approx\bs{0}}\leq g(\bs{y},\bs{z})|_{\bs{s}\approx\bs{0}}\leq u(\bs{y},\bs{z})|_{\bs{s}\approx\bs{0}}
	\end{align*}	
	for tolerable amounts of slip $\bs{s}\approx\bs{0}$. Replace $\bs{z}$ with a nominal value $\bs{Z}$ since the exact value of $\bs{z}$ is unknown to us. Lock the differentials if 
	\begin{align*}
	l(\bs{y},\bs{Z})\leq g(\bs{y},\bs{Z})\leq u(\bs{y},\bs{Z})
	\end{align*}
	does not hold. Unlock the differentials if a prespecified time $\Delta t$ has passed since locking and the vehicle has also moved a prespecified distance $\Delta d$.
\end{algo}


\subsection{Kinematic model of an articulated vehicle}
\label{sec:general}

\noindent The kinematics of load-haul-dump vehicles (a kind of low set articulated wheel loaders) is discussed in a literature on path-tracking in underground environments, see \emph{e.g.,} \cite{scheding1997slip,larsson1994navigating}. The model \cite{scheding1997slip} includes wheel slip angles. We cannot apply it directly to our setup since the vehicle geometry of wheel loaders and articulated haulers differ in key respects. A procedure for deriving kinematic models for the planar motion of articulated vehicles is presented in \cite{larsson1994navigating}. This paper generalizes the models of \cite{scheding1997slip,larsson1994navigating} by accounting for a broader range of ground vehicles, including but not limited to articulated haulers with wheel slip, using the technique from \cite{larsson1994navigating}.

\begin{figure}[htb!]
	\centering
	\includegraphics[width=0.5\textwidth]{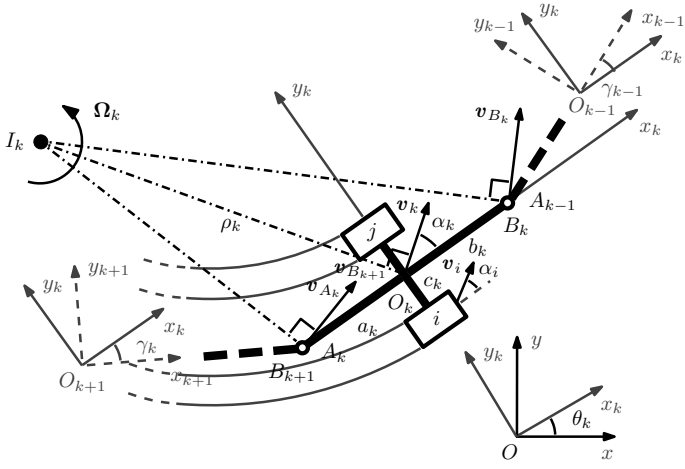}
	\caption[Kinematics of an articulated vehicle pair of wheels.]{Kinematics of an articulated vehicle. $I_k$ is the instantaneous center of rotation around which the $k$th axle rotates with velocity $\bs{\Omega}_k$. $\bs{v}_k$ is the translational velocity on the midpoint of the axle. The points $A_k$ and $B_k$ are on the longitudinal axle, situated at distances $a_k$ and $b_k$ from the origin $O_k$ of the local coordinate system. $\bs{v}_{A_k}$ and $\bs{v}_{B_k}$ are the point's translational velocities. Moreover, $c_k$ is half the axle track, $\rho_k=\Omega_k/v_k$ is the turning radius, $\alpha_k$ is the slip angle, and $\gamma_k$ is the angle between the coordinate axes $x_k$ and $x_{k-1}$. The angle $\theta_k$ separates the coordinate axis $x_k$ and the $x$-axis of a global (\ie interial fixed) coordinate system with origin $O$.}
	\label{fig:axle}
\end{figure}

Consider the $k$th transverse axle of an articulated vehicle with $n$ joints, as shown in Fig. \ref{fig:axle}. The velocity $\bs{v}_k=v_k\,[\cos\alpha_k\,\sin\alpha_k]^\mathsf{T}$ in the local coordinate origin $O_k$ is related to the velocities $\bs{v}_{A_k}$ and $\bs{v}_{B_k}$ in the points $A_k$ and $B_k$ by the equations $\bs{v}_{A_k}=\bs{v}_{k}+\bs{\Omega}_k\times\bs{r}_{O_kA_k}$ and $\bs{v}_{B_k}=\bs{v}_{k}+\bs{\Omega}_k\times\bs{r}_{O_kB_k}$. The velocity $\bs{v}_{B_{k+1}}$ is equal to $\bs{v}_{A_k}$ rotated by $\gamma_k$ degrees since $\gamma_k$ is the difference between the coordinates in $O_{k-1}$ and $O_k$, see Fig.  \ref{fig:axle}. These relations can be written on matrix form. First introduce a 2D screw, $\bs{\Psi}_k\equalperdefinition [v_k\,\,\Omega_k]^\mathsf{T}$, then
\begin{align*}
\bs{v}_{A_k}&=\bs{v}_{k}+\bs{\Omega}_k\times\bs{r}_{O_kA_k}\\
&=\begin{bmatrix}
v_k\cos\alpha_k\\
v_k\sin\alpha_k\end{bmatrix}+
\begin{bmatrix}1 &0 &0\\
0&1&0
\end{bmatrix}\begin{bmatrix}
0\\
0\\
\Omega_k\end{bmatrix}\times
\begin{bmatrix}
-a_k\\
0\\
0
\end{bmatrix}\\
&=\begin{bmatrix}
\cos \alpha_{k} & \matrixminusspace0\phantom{_k}\\
\sin \alpha_{k} & \matrixminus a_k
\end{bmatrix} \bs{\Psi}_k\equalperdefinition M_{A_k}\bs{\Psi}_k, \nonumber\\
\bs{v}_{B_k}&= \begin{bmatrix}
\cos \alpha_k & \matrixminusspace0\phantom{_k}\\
\sin \alpha_k & \matrixminusspace b_k
\end{bmatrix} \bs{\Psi}_k\equalperdefinition M_{B_k}\bs{\Psi}_k, \nonumber\\
\bs{v}_{B_{k+1}}&= \begin{bmatrix}
\cos\gamma_k & \matrixminus\sin\gamma_k\\
\sin\gamma_k & \matrixminusspace\cos\gamma_k
\end{bmatrix})\bs{v}_{A_k}\equalperdefinition R_{\gamma_k}\bs{v}_{A_k},\nonumber
\end{align*}
resulting in the system of equations
\begin{align}
M_{B_{k+1}}\bs{\Psi}_{k+1}=R_{\gamma_k}M_{A_k}\bs{\Psi}_{k}, \label{eq:kinematic1}
\end{align}
which may be solved for $\Omega_k$, $v_{k+1}$ and $\Omega_{k+1}$ as functions of $v_k$, $\gamma_k$, and $\dot{\gamma}_k$ by using the relation 
\begin{align}
\Omega_{k+1}=\Omega_k-\dot{\gamma}_k \label{eq:kinematic2}
\end{align}
and a few trigonometric identities. 

Moreover, the velocity $v_i$ of a single wheel with slip angle $\alpha_i$, expressed in the $k$th axle local coordinates, may be calculated as
\begin{align*}
\boldsymbol{v}_i&=\begin{bmatrix}
v_i\cos\alpha_i\\
v_i\sin\alpha_i
\end{bmatrix}=\bs{v}_k+\bs{\Omega}_k\times\bs{r}_{O_k i}\\
&=\begin{bmatrix}
v_k\cos\alpha_k\\
v_k\sin\alpha_k
\end{bmatrix}+\begin{bmatrix}
1 & 0 & 0\\
0 & 1 & 0
\end{bmatrix}\begin{bmatrix}
0\\
0\\
\Omega_k
\end{bmatrix}\times\begin{bmatrix}
0\\
-c_k\\
0
\end{bmatrix}\\
&=\begin{bmatrix}
\cos\alpha_k & c_k\\
\sin\alpha_k & 0
\end{bmatrix}
\boldsymbol{\Psi}_k 
\end{align*}
where $\alpha_i$ is the $i$th wheel slip angle. Note that the sign before $c_k$ implies that the wheel $i$ is on the right side of the vehicle.

\subsection{Kinematic model of an articulated hauler}
\label{sec:kinematic}

\begin{figure}[hbt!]
	\centering
	\includegraphics[width=0.4\textwidth]{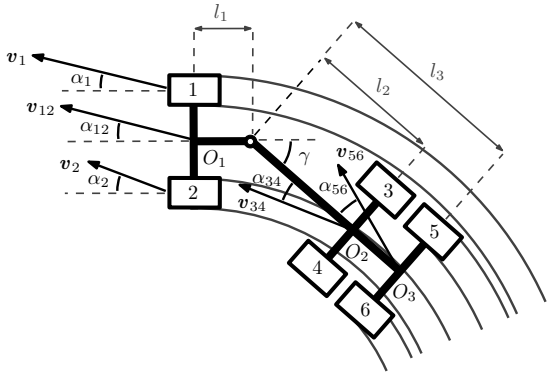}
	\caption[Variables and parameters of the kinematic vehicle model.]{	\label{fig:haulerBirdsview}Articulated hauler kinematics. The points $O_1$, $O_2$ and $O_3$ correspond to the locations of the mean wheels in a three-wheel articulated bicycle model. $1,\ldots,6$ are the actual wheels, $l_k$ is the distance from point $O_k,\,k\in\{1,2,3\}$ to the hinge, $\bs{v}_{ij}$ is the vehicle translational velocity at point $O_k$, $\gamma$ is the steering angle and $\alpha_{ij}$ is the $k$th mean wheel slip angle.} 
\end{figure}

\noindent Picture an articulated hauler as displayed in Fig.  \ref{fig:haulerBirdsview} and note the likeness with Fig. \ref{fig:axle}. The corresponding system of equations  \eqref{eq:kinematic1} and \eqref{eq:kinematic2} can be solved to yield
\begin{align}
v_{34}={}&\frac{(l_{2}\cos(\gamma+\alpha_{12})+l_{1}\cos\alpha_{12})v_{12}}{l_{2}\cos\alpha_{34}+l_{1}\cos(\gamma-\alpha_{34})} \nonumber \\
&+\frac{l_{1}l_{2}\sin\gamma\dot{\gamma}}{l_{2}\cos\alpha_{34}+l_{1}\cos(\gamma-\alpha_{34})} \label{eq:v34} \\
\nonumber\\
\Omega_1={}&\frac{\sin(\gamma+\alpha_{12}-\alpha_{34})v_{12}}{l_{2}\cos\alpha_{34}+l_{1}\cos(\gamma-\alpha_{34})}\nonumber\\
&+\frac{l_{2}\cos\alpha_{34}\dot{\gamma}}{l_{2}\cos\alpha_{34}+l_{1}\cos(\gamma-\alpha_{34})}
\label{eq:omega1}\\
\nonumber\\
\Omega_2={}&\frac{\sin(\gamma+\alpha_{12}-\alpha_{34})v_{12}}{l_{2}\cos\alpha_{34}+l_{1}\cos(\gamma-\alpha_{34})}\nonumber\\
&-\frac{l_{1}\cos(\gamma-\alpha_{34})\dot{\gamma}}{l_{2}\cos\alpha_{34}+l_{1}\cos(\gamma-\alpha_{34})},
\label{eq:omega2}
\end{align}
where the notation is that of Fig. \ref{fig:axle} and \ref{fig:haulerBirdsview}. 

Note that the equations \eqref{eq:v34}--\eqref{eq:omega2} encompasses the models of a wheel loader with and without slip in \cite{scheding1997slip} and \cite{larsson1994navigating} respectively as special cases.

Introduce the functions
\begin{align*}
p(\alpha_{12},\alpha_{34},\gamma)&=\frac{\sin(\gamma+\alpha_{12}-\alpha_{34})}{l_{1}\cos(\alpha_{34}-\gamma)+l_{2}\cos\alpha_{34}},\\
q(\alpha_{12},\alpha_{34},\gamma)&=\frac{l_{2}\cos(\alpha_{12}+\gamma)+l_{1}\cos\alpha_{12}}{l_{1}\cos(\alpha_{34}-\gamma)+l_{2}\cos\alpha_{34}},
\end{align*}%
where $p$ can be interpreted as the inverse of the tractor unit steady state turning radius (\emph{i.e.,} $\dot{\gamma}=0$) and $q$ as the quotient of the trailer and tractor unit steady state turning radii.

The equations \eqref{eq:v34}--\eqref{eq:omega2} may be tidied up by writing
\begin{align}
\hspace{-3mm}v_{34}&=q(\alpha_{12},\alpha_{34},\gamma)v_{12}+q(\matrixminus \pi/2,\alpha_{34},\gamma)\dot{\gamma}l_{1},\hspace{-1mm}\label{eq:v34q}\\
\hspace{-3mm}\Omega_1&=p(\alpha_{12},\alpha_{34},\gamma)v_{12}+p(\matrixminus\gamma+\pi/2,\alpha_{34},\gamma)\dot{\gamma}l_{2},\hspace{-3mm}\nonumber \\
\hspace{-3mm}\Omega_2&=p(\alpha_{12},\alpha_{34},\gamma)v_{12}+p(\matrixminus\pi/2,\alpha_{34},\gamma)\dot{\gamma}l_1,\nonumber 
\end{align}
with \eg $\dot{\gamma}l_{1}$ in equation \eqref{eq:v34q} being interpreted as a tractor unit velocity with a $\matrixminus\pi/2$ side slip angle.

\begin{figure}[htb!]
	\centering
	\includegraphics[width=0.45\textwidth, angle=-0.5]{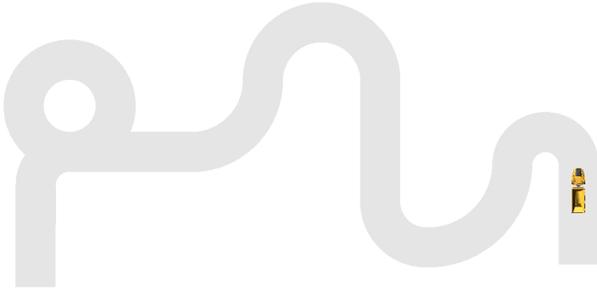}
	\caption{Test road in ADAMS. The articulated hauler runs from start to finish on the road centerline. The road lies in a plane (with the gravitational force $m\bs{g}$ as normal) and is mainly composed of curves.	The greatest circle sector curvature corresponds to the maximal steering angle of 45$^\circ$.}
	\label{fig:snake}
\end{figure}

To validate the kinematic model, we compare it to an existing articulated hauler model in MSC ADAMS, an environment for simulation of multibody dynamics. The ADAMS, model has previously been described in \cite{illerhag2000study,furmanik2015}. The articulated hauler in ADAMS is run repeatedly on a curved road, see Figure \ref{fig:snake}, while varying the gear and load mass. 

To validate the model \eqref{eq:v34q} we use two error metrics,
\begin{align}
d_{1}(x,\hat{x})&=\frac{1}{n}\sum_{k=1}^n{\left| x_k-\hat{x}_k\right|}, \label{eq:mean_error}\\
d_{\infty}(x,\hat{x})&=\max_{k\,\in\{1,\ldots,n\}}\left| x_k-\hat{x}_k \right|, \label{eq:max_error}
\end{align}
to measure the discrepancy between a true value $x$ and an estimate $\hat{x}$. The index $k$ is a discrete time instance and $n$ is the total number of time instances considered. Table \ref{tab:accuracy} displays the accuracy of an estimate
\begin{align*}
\hat{v}_{34}=q(\alpha_{12},\alpha_{34},\gamma)v_{12}+q(\nicefrac{\matrixminus \pi}{2},\alpha_{34},\gamma)\dot{\gamma}l_{1}
\end{align*}
compared to the velocity $v_{34}$ of the first bogie axle created as ouput from a simulation in  ADAMS. The three pairs of columns illustrate the error caused by setting $\dot{\gamma}$ and $\alpha_{ij}$ to zero. Observe that the error increases with gear choice (speed) and with load. Based on Table \ref{tab:accuracy} we conclude that only the models with non-zero steering angle velocity is the close to the true model. Slip angles are also helpful. The slip angles are not measured, although they can be estimated \cite{andersson2011estimation}. Hence we look for relations \eqref{eq:equation} that do not depend on slip angles.

\begin{table}[htb]
	\centering \small
	\captionsetup{singlelinecheck=off}
	\begin{tabular}{c@{\quad}c|c@{\quad}c|c@{\quad}c|c@{\quad}c}
		& & \multicolumn{2}{|l|}{$\alpha_{ij}\equalperdefinition0$,} & \multicolumn{2}{|l|}{$\alpha_{ij}\equalperdefinition0,$} & \multicolumn{2}{|l}{$\alpha_{ij}\in\mathbb{R},$}\\
		& & \multicolumn{2}{|l|}{$\,\,\,\,\,\dot{\gamma}\equalperdefinition 0$} & \multicolumn{2}{|l|}{$\,\,\,\,\,\dot{\gamma}\in\mathbb{R}$} & \multicolumn{2}{|l}{$\,\,\,\,\,\dot{\gamma}\in\mathbb{R}$}\\
		
		\hline
		\hline
		\phantom{${}^1$}\textbf{Load} & \textbf{Gear} & $d_1$ & $d_\infty$ & $d_1$ & $d_\infty$ & $d_1$ & $d_\infty$\\
		\hline	
		\hline
		
		\phantom{${}^1$}Zero	&	F1 & 1.9 & 20.9 & 0.4 & 2.1 & 0.1 & 0.6\\
		&	F2 & 2.9 & 31.7 & 0.6 & 3.0 & 0.2 & 0.7\\
		&	F3 & 3.7 & 41.3 & 0.7 & 3.5 & 0.2 & 0.7\\
		
		\phantom{${}^1$}Half	&	F1 & 2.0 & 21.0 & 0.4 & 2.3 & 0.1 & 1.0\\
		&	F2 & 3.0 & 32.3 & 0.6 & 3.2 & 0.2 & 1.4\\
		&	F3 & 3.8 & 40.6 & 0.7 & 3.9 & 0.2 & 1.6\\
		
		\phantom{${}^1$}Full	&	F1 & 2.1 & 22.2   & 0.4 & 2.6 & 0.3 & 1.7\\
		&	F2 & 3.1 & 32.9   & 0.6 & 3.5 & 0.3 & 2.4\\
		&	F3 & 3.9 & 42.1   & 0.8 & 4.1 & 0.3 & 2.5
	\end{tabular}
	\caption[Mean and maximal relative trailer unit velocity estimation errors.]{Mean and maximal relative trailer unit speed estimation errors, equations \eqref{eq:mean_error} and \eqref{eq:max_error} with $x=v_{34}$ (as measured in ADAMS) and $\hat{x}$ given by equation \eqref{eq:v34q}. The units are cm s$^{-1}$. The road of Fig. \ref{fig:snake} is used. Time instances with steering angles less than $1^\circ$ were removed from consideration as they correspond to straight parts of the test road. The rows and columns vary: 
		\begin{itemize}
			\item[(i)] setting the slip angles $\alpha_{12},\alpha_{34}$ to zero,
			\item[(ii)] the three first forward gears F1, F2 and F3 (roughly corresponding to $v_{12}=2,\,2.5$ and $3$ m s$^{-1}$ respectively)
			\item[(iii)] zero, half, or full load (39 000 kg),
			\item[(iv)]  using a steady state turning model ($\dot{\gamma}=0$) or a transient state turning model ($\dot{\gamma}\in \mathbb{R}$).
		\end{itemize}	
		\label{tab:accuracy}}
\end{table}

\subsection{Basic sensor network}
\label{sec:basic}


\noindent The basic sensor network is based on \cite{atc} and consists of a steering angle sensor for the angle between the tractor and the trailer unit and four tachometers that measure rotational velocities at various points on the driveline (see Section \ref{sec:hauler} and Fig. \ref{fig:atc}). There are no tachometers on the wheels since they engage directly with the terrain, which would drastically affect the expected lifespan of any sensors placed there. 

In practice, these five sensors provide enough information for a complete traction control algorithm \cite{atc}. However, we will also consider algorithms based on additional information. Of the basic sensors, the steering angle $\gamma$ (and its derivative $\dot{\gamma}$) is the most interesting since its readings relate to our kinematic model of the hauler, equation \eqref{eq:v34}. Traction control for the case of $\gamma=0$ is fairly straightforward from a theoretical point of view  (see \cite{furmanik2015} for details) so we focus on the case of $\gamma\neq0$.

Recall that our traction control approach described in Section \ref{sec:control} is based on comparing quanties that are affected and not affected by slip. The kinematic equation \eqref{eq:v34q} can be used to obtain
\begin{align}
g(\bs{y},\bs{z})|_{\bs{z}=\bs{0}}={}&\frac{\omega_3+\omega_4}{2}r-q(0,0,\gamma)\frac{\omega_1+\omega_2}{2}r-\nonumber\\
&q(\pi/2,0,\gamma)l_1\dot{\gamma}\nonumber\\
={}&\omega_{\textit{bg, in}}r/i-q(0,0,\gamma)\omega_{\textit{dbx, out}}r/i-\nonumber\\
&q(\pi/2,0,\gamma)l_1\dot{\gamma}\label{eq:basic}
\end{align}
where $i$ is the gear conversion ratio from the differentials over the hub reductions to the wheels ($i=i_{\textit{diff}}\cdot i_{\textit{hub}}=3.09\cdot6=18.54$), $\bs{y}=[\omega_{12}\,\omega_{34}\,\gamma\,\dot{\gamma}]^\textsf{T}$, $\bs{z}=[\alpha_{12}\,\alpha_{34}]^\textsf{T}$. Note that we used $(\omega_i+\omega_j)r/2=v_{ij}$ as we have set $\bs{z}=\bs{0}$.

The error in equation \eqref{eq:basic} should not be larger than those in Table \ref{tab:accuracy} with the slip angles $\alpha_{12}$ and $\alpha_{34}$ set to zero. If they are we can conclude that it is due to wheel slip and engage the differential locks. Note that the errors for setting $\dot{\gamma}=0$ in Table \ref{tab:accuracy} are rather large. As such, with only basic sensors, it is important to use an accurate kinematic model. This explains why the VCE ATC system includes a steering angle sensor.

\subsection{Ground speed sensor}
\label{sec:ground}

\noindent GPS receivers and ground speed radars are examples of sensors that can be used to estimate or measure the ground speed. With a two-antenna GPS reciever, one antenna on the tractor unit and one on the trailer, it is possible to estimate wheel slip angles \cite{andersson2011estimation}. Note that the 1 Hz sample rate of a GPS is rather slow, but  the ground speed can be assumed piece-wise constant. What is worse, a hauler may be put to work in GPS denied enviroments, \eg mines or tunnels. Still, a GPS is probably a more likely option than a ground speed radar due to the former's multi-purpose versatility. A ground speed radar also needs a clear line of sight to the ground. It is hence exposed to mud and can be rendered inoperable.

Output from a ground speed sensor can be compared to tachometer readings of the driveline shaft's angular speed to calculate the slip according to definition \eqref{eq:deltaomega}. An algorithm that knows the vehicle ground speed would hence be able to detect most occurences of wheel slip.


Assume that the speed is measured somewhere on the tractor unit and that it is recalculated to correspond to the mean 1st and 2nd wheel speed (see Fig. \ref{fig:haulerBirdsview}). The expression
\begin{align}
\frac{\omega_{1}+\omega_{2}}{2} r-v_{12}=\omega_{\textit{dbx, out}}r/i-v_{12}, \label{eq:v12_difference}
\end{align}
where $i$ is the gear conversion ratio from the dropbox output to the wheel ($i=i_{\textit{diff}}\cdot i_{\textit{hub}}=3.09\cdot6=18.54$ ), is an equation of the type \eqref{eq:equation} which can be rewritten as a system of inequalities of the type \eqref{eq:condition} to detect a slip of the front mean wheel. Note that we have set  $\bs{z}=[\alpha_{12}\,\alpha_{34}]^\mathsf{T}=\bs{0}$. A slip of one of the front bogie wheels is detected as
\begin{align}
\smash{\frac{\omega_{3}+\omega_{4}}{2}} r-v_{34}={} &\omega_{\textit{bg,\,in}}r/i-q(0,0,\gamma)v_{12}+\nonumber\\
&q(\nicefrac{\matrixminus \pi}{2},0,\gamma)\dot{\gamma}l_{1},\label{eq:v34_difference}
\end{align}
where $i$ is a gear conversion ratio and $\omega_{bg,in}$ is the input tachometer on the bogie axle, see Fig. \ref{fig:atc}.



\subsection{Individual wheel tachometers}

\noindent The VCE ATC system tachometers (see Fig. \ref{fig:atc}) do not measure wheel angular speeds, but rather the angular speed of driveline shafts. These speeds are proportional to the speeds of the mean wheels under the assumption of a rigid driveline. Measurements of the wheel angular speeds facilitate slip detection but are difficult to carry out since tachometers placed at those locations have a short expected lifetime. One mud bath is enough to kill the sensors. To protect them, they need to be encapsulated into the wheel hub. In this section we assume that output from tachometers in all wheel hubs is readily available.

Consider wheel $i$ in Fig. \ref{fig:axle}. From $v_i=\Omega_k\rho_i$ we get $v_i\cos\alpha_i=\Omega_k\rho_i\cos\alpha_i$. Moreover, by definition \eqref{eq:deltaomega}, $s_{i}=\omega_i r- v_i\cos\alpha_i=\omega_ir-\Omega_k\rho_i\cos\alpha_i$. Calculate
\begin{align*}
s_{1}-s_{2}&=\omega_1 r-\omega_2r-\Omega_k(\rho_1\cos\alpha_1-\rho_2\cos\alpha_2)\\
&=(\omega_1 -\omega_2)r+2c_1\Omega_k,
\end{align*}
where we used the geometry in Fig. \ref{fig:axle} to find  $\rho_1\cos\alpha_1-\rho_2\cos\alpha_2=-2c_1$. A similar equation holds for the trailer,
\begin{align*}
s_{3}-s_{4}&=(\omega_3-\omega_4)r+2c_3\Omega_{k+1}.
\end{align*}
Combine the equations for the tractor and trailer unit to obtain
\begin{align}
s_{1}-s_{2}-s_{3}+s_{4}=(\omega_1-\omega_2-\omega_3+\omega_4)r+2c\dot{\gamma} \label{eq:1234}
\end{align}
where we used $\Omega_2=\Omega_1-\dot{\gamma}$ and $c_1=c_3=c$ \cite{spec}. Note that we are able to remove the wheel slip angles $\alpha_{i}$ without setting them to zero as in Section \ref{sec:basic} and \ref{sec:ground}. The equation \eqref{eq:1234} is precisely the type of expression \eqref{eq:equation} on the form $g(\bs{y},\bs{z})|_{\bs{s=0}}=g(\bs{y})_{\bs{s=0}}$ that we want.

To find all combinations of wheel slip that cannot be detected we write equation \eqref{eq:1234} on matrix form, 
\begin{align*}
\begin{bmatrix}
1 &\textit{--}1 & \textit{--}1 & 1
\end{bmatrix}
\begin{bmatrix}
s_1\\
s_2\\
s_3\\
s_4
\end{bmatrix}=0
\end{align*}
and calculate the null space as
\begin{align}\label{eq:ker}
\mathrm{ker}\begin{bmatrix}
1 &\textit{--}1 & \textit{--}1 & 1
\end{bmatrix}= \linspan \left\{
\begin{aligned}
\left[
\begin{aligned}
1\\
1\\
0\\
0
\end{aligned}
\right],
\left[
\begin{aligned}
0\\
0\\
1\\
1
\end{aligned}
\right],
\left[
\begin{aligned}
\textit{--}1\\
1\\
\textit{--}1\\
1
\end{aligned}
\right]
\end{aligned}
\right\}.
\end{align}
Slip is undetectable to the control strategy whenever
\begin{align*}
\bs{s}=[s_1\,s_2\,s_3\,s_4]^\mathsf{T}\in\mathrm{ker}\begin{bmatrix}
1 &\textit{--}1 & \textit{--}1 & 1
\end{bmatrix}.
\end{align*}
Note that the nullspace \eqref{eq:ker} is the span of vectors where multiple wheels slip. Hence \eqref{eq:1234}, like our other criteria, is the most useful for detecting the slip of a single wheel.

\subsection{Tuning}
\label{sec:tuning}

\noindent The functions $u(\bs{y})$ and $l(\bs{y})$ in equation \eqref{eq:condition} 
should be tuned so that no slip is mistakenly detected during normal driving (\ie no false positives) without delaying the detection of actual slip (\ie minimize the amount of false negatives). This tuning is done based on data from the simulation model in ADAMS. 

\begin{figure}[htb!]
	\centering
	\includegraphics[width=0.45\textwidth, angle=-0.5]{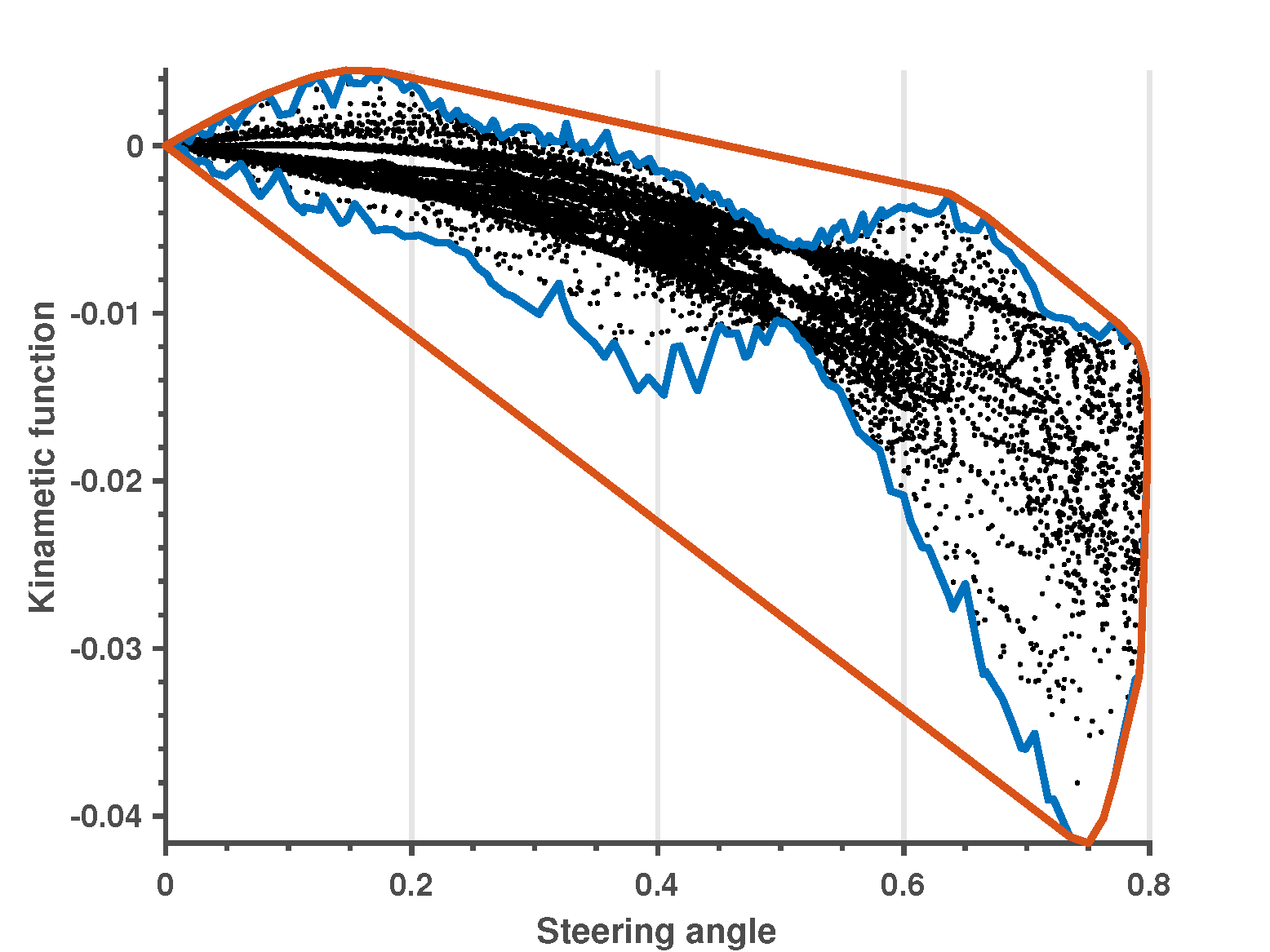}
	\caption{Pairs $(\gamma,g(\bs{y},\bs{z})|_{\bs{z}=\bs{0}})$, where $g(\bs{y},\bs{z})|_{\bs{z}=\bs{0}}$ is given by \eqref{eq:basic}, their boundary (blue), and convex hull (red). The data is obtained from the ADAMS model on the road in Fig. \ref{fig:snake} on the gears 1,  2, and 3 with empty, half, and full load.}
	\label{fig:basic_boundary}
\end{figure}

Consider first the basic sensor network. Figure \ref{fig:basic_boundary} displays a point cloud of pairs $(\gamma,g(\bs{y},\bs{z})|_{\bs{z}=\bs{0}})$ where $g(\bs{y},\bs{z})|_{\bs{z}=\bs{0}}$ is given by the kinematic relation between the tractor and trailer unit \eqref{eq:basic}. The road in Fig. \ref{fig:snake} on which the data is collective is such that only tolerable amounts of wheel slip are present. The upper boundary of the point cloud is $u(\gamma)$, the lower boundary is $l(\gamma)$. These functions could also be assigned to upper and lower curves of the convex hull of the point cloud, but that seems overly conservative, at least for the lower function $l(\gamma)$. Both functions are piece-wise linear. As is clear by inspection, the relation is not symmetric, \emph{i.e.,} $l(\gamma)\neq-u(\gamma)$. It is possible to  approximate the bounds using two quadratic functions $u(\gamma)=u_1\gamma+u_2\gamma^2$ and $l(\gamma)=l_1\gamma+l_2\gamma^2$ although a piece-wise linear map based on stored data is also an option. 

\begin{figure}[htb!]
	\centering
	\includegraphics[width=0.45\textwidth, angle=-0.5]{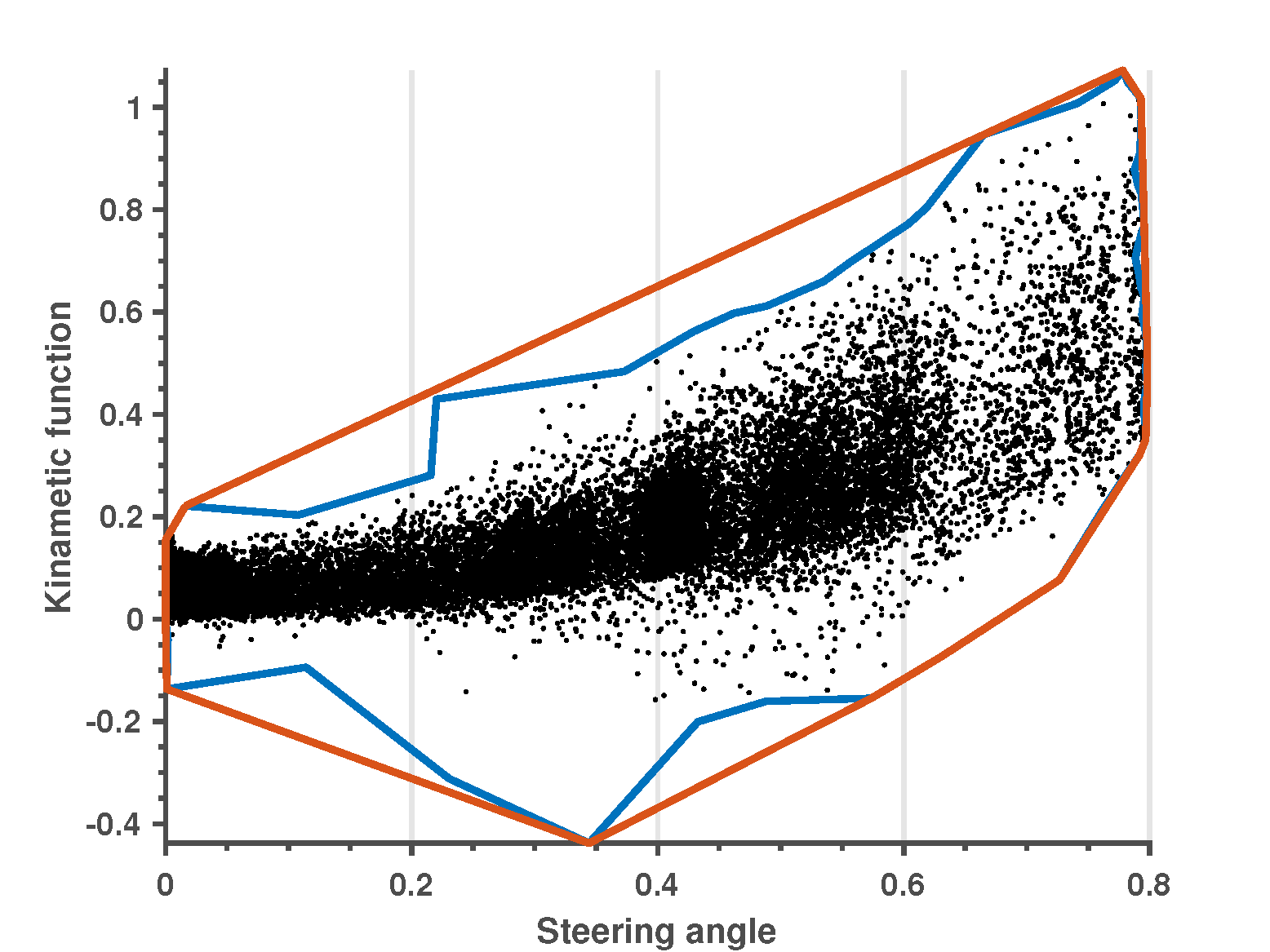}
	\caption{Pairs $(\gamma,g(\bs{y},\bs{z})|_{\bs{z}=\bs{0}})$ , where $g(\bs{y},\bs{z})|_{\bs{z}=\bs{0}}$ is given by \eqref{eq:v12_difference}, their boundary (blue), and convex hull (red). The data is obtained from the ADAMS model on the road in Fig. \ref{fig:snake} on the gears 1,  2, and 3 with empty, half, and full load.}
	\label{fig:ground_boundary}
\end{figure}

Consider the ground speed sensor algorithm, where $g(\bs{y},\bs{z})$ is given by \eqref{eq:v12_difference} or \eqref{eq:v34_difference}. We only consider equation \eqref{eq:v12_difference} here since the case of \eqref{eq:v34_difference} is similar. Figure \ref{fig:ground_boundary} shows that the relation between $\gamma$ and \eqref{eq:v12_difference} is that of a quadratic function. Hence it is possible to fit two functions $u(\gamma)=u_0+u_1\gamma+u_2\gamma^2$ and $l(\gamma)=u_0+u_1\gamma+u_2\gamma^2$ to describe the upper and lower bounds. Unlike for the basic sensor network based algorithm, here we need to use $u(0)=u_0\approx0.2$ and $l(0)=l_0\approx-0.2$. It may be possible to lower the $l_0$ value by using a filter since the value of $l(0)=l_0\approx-0.2$ appears to be an outlier. Still, $l_0$ is unlikely to be a problem since it would indicate that the vehicle is braking due to $v>\omega r_0$

\begin{figure}[htb!]
	\centering
	\includegraphics[width=0.45\textwidth, angle=-0.5]{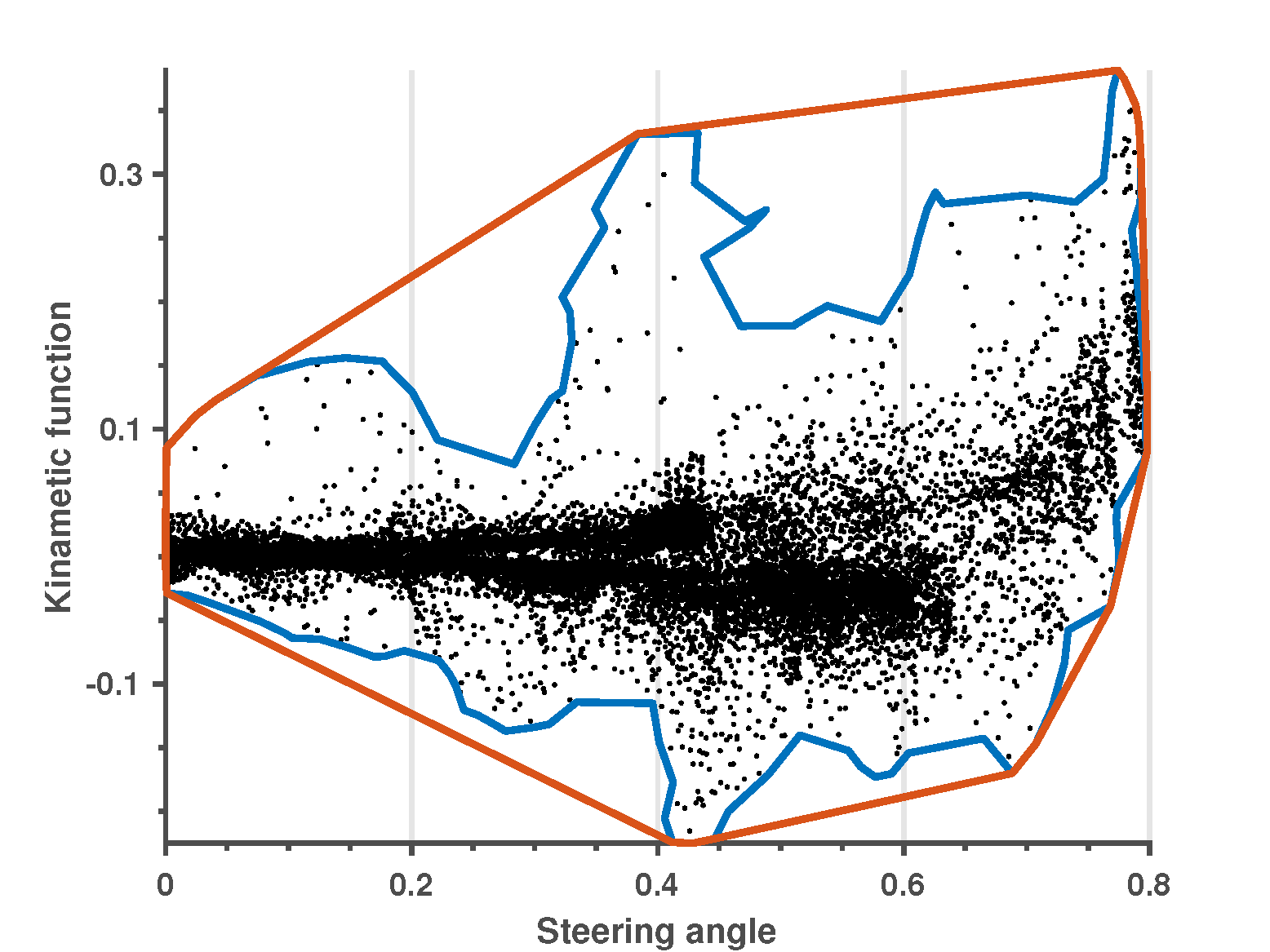}
	\caption{Pairs $(\gamma,g(\bs{y},\bs{z})|_{\bs{z}=\bs{0}})$ , where $g(\bs{y},\bs{z})|_{\bs{z}=\bs{0}}$ is given by \eqref{eq:1234}, their boundary (blue), and convex hull (orange). The data is obtained from the ADAMS model on the road in Fig. \ref{fig:snake} on the gears 1,  2, and 3 with empty, half, and full load.}
	\label{fig:wheels_boundary}
\end{figure}

Consider the individual wheel tachometers algorithm, where $g(\bs{y})$ is given by \eqref{eq:1234}. The $\dot{\gamma}$ signal in ADAMS changes rapidly. To remove some sharp peaks from the curve we run a filter that averages the kinematic relation to equal the mean of its value over the last five points. Still, the shape of $l(\gamma)$ and $u(\gamma)$ given by Fig. \ref{fig:wheels_boundary} is irregular. In theory, the condition \eqref{eq:1234} is very attractive since it excludes the slip angles. In practice we see that $l(\bs{y})$ and $u(\bs{y})$ has to be chosen quite large. In particular, we have $l(\bs{0}),u(\bs{0})\neq 0$ which was not the case for the relation \eqref{eq:basic} based on the basic sensor network. Moreover, the boundary is very irregular in some places, suggesting that more data is required. Suppose the error $g(\bs{y})$ is equal to the slip of wheel $1$. Then we have $s_1=g(\bs{y})$. The RPM difference over the front transversal differential is $6g(\bs{y})$ where the $6$ stems from the hub reduction gear conversion ratio. We have to tolerate a maximum diffrence of $u(\pi/4)=6\cdot0.35*60/(2\pi r)\approx 21$ RPM. The tolerable difference for $\gamma=0$ is at least $u(0)=6\cdot0.1*60/(2\pi r)\approx 6$ RPM, $l(0)=6\cdot0.05*60/(2\pi r)\approx 3$ RPM.
 
\section{Conclusions}

\noindent Three algorithms are presented together with their relative strengths and weaknesses. The basic sensor network suffices for a feasible automatic traction control (ATC) algorithm, seeing as one is implemented on Volvo CE haulers \cite{atc}. It may be possible to improve on this algorithm using additional sensor data. In particular, ground speed measurements based on GPS reciever data is expected to be available in the future. We can also speculate on the benefits of individual wheel tachometers. In theory, individual wheel tachometers present the most accurate data, allowing us to remove the influence on wheel slip angles. In practice however, we find that the tachometer data is plagued by too much variation and many outliers. The ground speed sensor data is well-behaved but offers no significant advantage over the basic sensor network based algorithm. In conclusion, the results of this study favor the basic sensor network based algorithm. It should be noted that other algorithms besides the ones presented here could be developed, and, broadly speaking, it is preferable to have more sensor output as compared to less.

\section{Acknowledgement}

\noindent This paper is based on the author's master's thesis carried out at Volvo Construction Equipment (VCE) in Eskilstuna, Sweden during 2009--2010. The work was jointly funded by VCE and the Swedish innovation agency Vinnova. Thanks to Gianantonio Bortolin, Ulf Andersson, and Ulf J\"onsson.

{\footnotesize\bibliography{autosam}}          



\end{document}